\documentclass[12pt, a4paper]{article}
\usepackage[left=1in,right=1in,top=1in,bottom=1in]{geometry}
\usepackage[english]{babel}
\usepackage{parskip}
\usepackage[utf8]{inputenc}
\usepackage[T1]{fontenc}
\usepackage{mathptmx}
\usepackage{graphicx}
\usepackage{titlesec}
\usepackage{pdfpages}
\usepackage[version=4]{mhchem}
\usepackage{amsfonts}
\usepackage[font=small]{caption}
\usepackage{amsmath}
\usepackage{mathtools}
\usepackage{nicefrac}
\usepackage[thinc]{esdiff}
\usepackage{derivative}
\usepackage{siunitx}
\usepackage{array}
\usepackage{booktabs}
\usepackage{multirow}
\usepackage[affil-it]{authblk}
\usepackage[colorinlistoftodos]{todonotes}
\usepackage{hyperref}
\usepackage{lineno}	 
\usepackage{csquotes}

\usepackage[doi=false,isbn=false,url=false,eprint=false,backend=biber,style=numeric,giveninits=true,sorting=none,maxnames=50]{biblatex}

\usepackage{setspace}
\usepackage{rotating}					
\usepackage{geometry}					
\usepackage{longtable}					
\usepackage{parcolumns}					
\usepackage{tabularx} 					
\usepackage{threeparttable}
\usepackage{rotating}
\usepackage{caption}
\usepackage{subcaption}
\usepackage{amsthm}
\usepackage[official]{eurosym}
\DeclareSIUnit{\euro}{\mbox{\euro}}

\usepackage{xcolor}
\definecolor{nile}{RGB}{26,135,133}
\DeclareRobustCommand\colorSS {\textcolor{nile}}
\definecolor{blue}{RGB}{2,133,204}
\DeclareRobustCommand\colorDistAll {\textcolor{blue}}
\definecolor{blue3}{RGB}{49,133,155}
\DeclareRobustCommand\colorDist {\textcolor{blue3}}
\definecolor{brown}{RGB}{212, 120, 0}
\DeclareRobustCommand\colorCentFinal {\textcolor{brown}}
\definecolor{purple}{RGB}{83,3,99}
\DeclareRobustCommand\colorCentInter {\textcolor{purple}}
\definecolor{black}{RGB}{0,0,0}

\newcommand{\CentSS}{OptSS}
\newcommand{\CentFinal}{FOwFinalTank-CO}
\newcommand{\CentInter}{FOwInterTank-CO}
\newcommand{\DistAll}{FOwInterTank-DO}
\newcommand{\DistRSR}{FOforRSR-DO}
\newcommand{\DistB}{FOforB-DO}
\newcommand{\DistG}{FOforG-DO}

\makeatletter
\newcommand*{\rom}[1]{\expandafter\@slowromancap\romannumeral #1@}
\makeatother

\titleformat{\section}
{\large\bfseries}
{\thesection.}{0.5em}{}
\titlespacing{\section}{0mm}{1.3pc}{0.8pc}

\titleformat{\subsection}
{\normalfont\bfseries}
{\thesubsection.}{0.5em}{}
\titlespacing{\subsection}{0cm}{0.8pc}{0.4pc}

\titleformat{\subsubsection}
{\normalfont\bfseries}
{\thesubsubsection.}{0.5em}{}
\titlespacing{\subsubsection}{0cm}{0.8pc}{0.4pc}

\titleformat*{\paragraph}{\normalfont\bfseries}
\titlespacing{\paragraph}{0cm}{0.5pc}{0.8pc}

\titleformat*{\subparagraph}{\normalfont\bfseries}
\titlespacing{\subparagraph}{0cm}{0pc}{0.8pc}

\providecommand{\keywords}[1]
{
	\small	
	\textbf{\textit{Keywords---}} #1
}

\makeatletter
\renewcommand{\maketitle}{\bgroup\setlength{\parindent}{0pt}
	\begin{flushleft}
		\textbf{\@title}
		
		\@author
	\end{flushleft}\egroup
}
\makeatother

\title{Optimal Design and Flexible Operation of a Fully Electrified Biodiesel Production Process}
\author[1]{\small Mohammad~El~Wajeh}
\author[1]{Adel~Mhamdi}
\author[2,1,3,*]{Alexander~Mitsos}

\affil[1]{Process Systems Engineering (AVT.SVT), RWTH Aachen University, 52074 Aachen, Germany}
\affil[2]{JARA-CSD, 52056 Aachen, Germany}
\affil[3]{Energy Systems Engineering (IEK-10), Forschungszentrum Jülich, 52425 Jülich, Germany}

\affil[*]{Corresponding author: \href{mailto:amitsos@alum.mit.edu}{amitsos@alum.mit.edu}}

\date{}

\begin{document}
	\onehalfspacing
	{\fontsize{16pt}{18pt}\selectfont \maketitle}	
	\vspace{0.3cm}
	\vspace{0.8cm}
	\begin{abstract}
	\vspace{-0.4cm}
	\noindent		
	The flexible operation of electrified chemical processes, powered by renewable electricity, offers economic and ecological incentives, paving the way for a more sustainable chemical industry. However, this requires a departure from the traditional steady-state operation, posing a challenge to process design and operation. We propose an approach for flexible biodiesel production by utilizing intermediate and final buffer tanks. These tanks serve to decouple the dynamics between different process unit operations. We investigate the flexibility of three different process configurations and compare the outcomes of their dynamic optimization strategies with those of a steady-state operation, considering a typical demand-side response scenario. For all strategies, we employ local gradient-based optimization and solve them using our open-source optimization framework DyOS. Our findings indicate that by incorporating intermediate buffer tanks, we fully exploit the process flexibility potential, leading to total energy savings of up to \SI{29}{\percent}. We also propose an implementation based on distributed optimization problems, yielding similar savings while significantly reducing computational costs by over 10-fold and improving convergence, thus paving the way toward online application. Overall, this study highlights the potential benefits of incorporating buffer tanks and optimization strategies in the flexible operation of electrified chemical processes, even those initially designed solely for steady-state operation.		
	\end{abstract}

	\newpage

	\section{Introduction} \label{sec:introduction}

    A transition toward an electrified, renewable-powered, and flexibly-operated chemical industry contributes to sustainable chemical production \cite{Schiffer.2017,Mitsos.2018,Barton.2020}. Flexible chemical processes can be operated dynamically, allowing for continuous adjustments in production rates to optimize energy usage and take advantage of fluctuating energy prices and availability \cite{Zhang.2016,Agora}. Flexibility can also apply to product purities and types, with varying products and grades requiring different power levels, and diverse feedstocks available \cite{Bruns.2020,Cegla.2023}. However, most conventional chemical processes are designed for continuous operation around steady-state (SS) operating points, which necessitates a constant energy supply, thereby posing significant challenges for flexible operation \cite{Cegla.2023}. Therefore, this paradigm shift requires a reevaluation of the conventional process design and operation to enable greater flexibility in accommodating variable energy availability.
	
	Large-scale industrial processes operate at time scales that are typically comparable to the frequency of energy price fluctuations. Consequently, the plant may remain in a transient state for prolonged periods during dynamic operation. Thus, an optimal flexible operation must consider process efficiencies, operational limits, and product qualities. This necessitates trajectory optimization, which entails optimizing the process degrees of freedom over a specified time horizon to maximize profit or minimize carbon footprint while satisfying operational constraints and meeting product quality requirements \cite{Cegla.2023}. Several dynamic optimization (DO) techniques, including direct sequential \cite{Sargent.1978,Kirches.2012} and full discretization methods \cite{Kameswaran.2006}, enable such trajectory optimization. The potential of optimal flexible trajectories has been investigated in several load-shifting applications for electricity-intensive processes like air separation units \cite{Caspari.2019c,Caspari.2020b}, water desalination plants \cite{Ghobeity.2010,K.Oikonomou.2020}, chlor-alkali electrolyzers \cite{Bree.2019,Otashu.2019}, and multi-energy systems \cite{Pablos.2021,Baader.2022}. However, numerical solutions to these optimization problems can be challenging, with computational costs increasing with model size and complexity, length of considered time horizon, and temporal resolution.
		
	Chemical plants are composed of several unit operations, most prominently involving reaction and separation. These different unit operations exhibit varying levels of operational flexibility, making it difficult to utilize the plant's full flexibility potential or adjust its load as a whole. For instance, the production flexibility of a unit operation may be limited by the operational constraints, particularly level limits, of the downstream processes. Incorporating intermediate and final storage units for the products may allow individual units to operate at different flexibility levels within their respective limits. Such a solution is simple yet effective as it also avoids the need to retrofit the unit operation sizes. Nonetheless, the exchange between upstream and downstream processing through these storage units must be coordinated optimally. This operational complexity and flexibility-oriented process design have not been fully comprehended in existing chemical processes involving reaction, separation, and recycle parts that consider some degree of process flexibility \cite{Hank.2018,Chen.2021}. Despite the added complexity, this strategy may become one of the new paradigms for chemical process design and operation in the era of renewable-powered chemical production.
	
	Biodiesel production processes are an exemplary case of classical chemical plants that involve reaction and separation processes along with multiple material recycle streams \cite{Zhang.2003,ElWajeh.2023}. Moreover, biodiesel is a ``renewable'' fuel derived from biomass that has the potential to replace fossil-based diesel \cite{Zhang.2003}. However, its production costs are higher than those associated with producing conventional diesel \cite{Gebremariam.2018}. Thus, reducing total manufacturing costs and leveraging fluctuating energy prices through flexible operations are critical to ensuring its economic viability. Therefore, designing an optimal and flexible electrified biodiesel production process would not only explore demand-side management in a classical chemical plant but also unify sustainability concepts from various research areas. This integration is a key aspect for future industrial biofuel and chemical synthesis using renewable energy.

    We investigate a fully electrified process for the homogeneous transesterification of vegetable oil using an alkali catalyst to produce biodiesel. Our process design builds upon our recently published work \cite{ElWajeh.2023}, in which we developed a rigorous first-principle model in Modelica with two plantwide control structures, and made it available as open-source. The process has a main reaction part, several separation processes, and two recycle streams. It produces two final products with specific quality requirements, as glycerol is a by-product. The process involves operational limits and stringent quality standards, which pose challenges to exploiting process flexibility. As the process educts are in the liquid state, they are suitable for intermediate and final storage without requiring additional liquefaction. Moreover, the process provides an opportunity to explore how operational flexibility affects heat integration across multiple units.

    To investigate the potential of flexibility-oriented designs and operations of the process, we examine three different process configurations using intermediate and final buffer tanks and compare the operating profits obtained from their respective economic DO strategies with that of a conventional optimal SS operation. The three process configurations differ based on the number and the location of the buffer tanks incorporated. We employ the three DO strategies offline, utilizing local gradient-based optimization. Furthermore, we demonstrate that intermediate buffer tanks not only decouple dynamics between different process parts but also facilitate the implementation of a distributed optimization strategy with smaller problem sizes, leading to enhanced computational performance.
	
	The manuscript is structured as follows. First, we introduce the biodiesel production process under consideration and provide a summary of our modeling approach. Next, we discuss the several configurations examined to achieve process flexibilization, along with the corresponding optimization problem assumptions and formulations. We further elaborate on the considered operational scenario and implementation before presenting and discussing the results. Lastly, we draw conclusions based on our findings.	
	
	\section{Biodiesel and glycerol production process} \label{sec:process}

	Figure~\ref{fig:PFD} depicts the flowsheet of the entire process, including all buffer tanks and optimization variables. In Section~\ref{sec:optimization}, we elaborate on the use of buffer tanks and optimization strategies. For completeness, we present a summary of the process description in this section, excluding buffer tanks, and briefly discuss our modeling methodology and underlying assumptions. We refer the reader to our prior work \cite{ElWajeh.2023} for more detailed information on the process description and modeling.
	
	\subsection{Process flowsheet} \label{sec:flowsheet}

    The oil feed composition consists of $95$~wt\unit{\percent} triolein and $5$~wt\unit{\percent} diolein \cite{Zhang.2003,ElWajeh.2023}. We use methanol for transesterification and sodium hydroxide (\ce{NaOH}) solution as a catalyst. They are fed into the transesterifier, which converts oil and methanol into fatty acid methyl ester (FAME) and glycerol. The resulting products are then separated and purified. The methanol column recovers the unreacted methanol to recycle it to the reactor. The bottom product is cooled and sent to a decanter to separate most of the glycerol, water, the dissolved \ce{NaOH}, and methanol from FAME and unreacted oil. The light product is then further purified in the FAME column, producing biodiesel that meets European quality standards \cite{BritishStandardsInstitution.2010}. The residual oil is recycled back into the reactor after being mixed with fresh oil. The decanter bottom is sent into a neutralizer unit, where it is neutralized using a phosphoric acid solution to remove dissolved \ce{NaOH} species \cite{Myint.2009}. The solution is then filtered to remove the formed salt in a filter unit. The remaining liquid is purified in a glycerol column to remove water and methanol from glycerol with a pharmaceutical-grade purity of $99$~wt\unit{\percent} \cite{Kariwala.2012, ElWajeh.2023}.

	\subsection{Process modeling} \label{sec:model}
	
	Before presenting the electrification approach, we briefly summarize the considered modeling formulation and assumptions for the unit operations and thermodynamics. In addition, we provide the base-layer control considered and further modeling aspects.

	\subsubsection{Unit operations}
	
	We describe the transesterifier as a perfectly-mixed CSTR using energy and material balances. We model the reaction system based on the well-studied reversible three-step transesterification system in the literature \cite{Noureddini.1997,ElWajeh.2023}. We model the transesterifier jacket as a series of equivalent CSTRs. We utilize MESH models for the decanter and each stage of the distillation columns. Our approach includes several assumptions, such as the perfect mixing of vapor and liquid phases, and that the tray holdup only accounts for the liquid phase since the vapor holdup is typically insignificant \cite{Raghunathan.2004}. We use a quasi-steady-state approximation for the energy balances at the trays, which allows reformulating the model to an index-1 differential-algebraic equation (DAE) system \cite{Raghunathan.2004}.
	
	\subsubsection{Chemical system and thermodynamics}
	
	We consider ten species, including triolein, diolein, methanol, water, \ce{NaOH}, and phosphoric acid, as feeds. Monoolein is an intermediate educt. The products are methyl oleate (\ce{FAME}), glycerol, and monosodium phosphate (solids product). We use the DIPPR temperature-dependent correlations to determine the molar heat capacities \cite{DesignInstituteforPhysicalPropertyData.2010}. We then determine species enthalpies, entropies, and Gibbs free energies through analytical integration. We determine the vapor pressures using the extended Antoine correlation \cite{Poling.2001}. To account for non-ideality in the liquid phase, we use the NRTL model \cite{Renon.1968}. The Racket equation is used to determine liquid-mixture molar densities \cite{Rackett.1970}. We employ isofugacity conditions to describe the liquid-liquid and vapor-liquid equilibrium, assuming that the number of existing phases is known.

    \begin{figure}[t]
		\centering
		\includegraphics[width=1\linewidth]{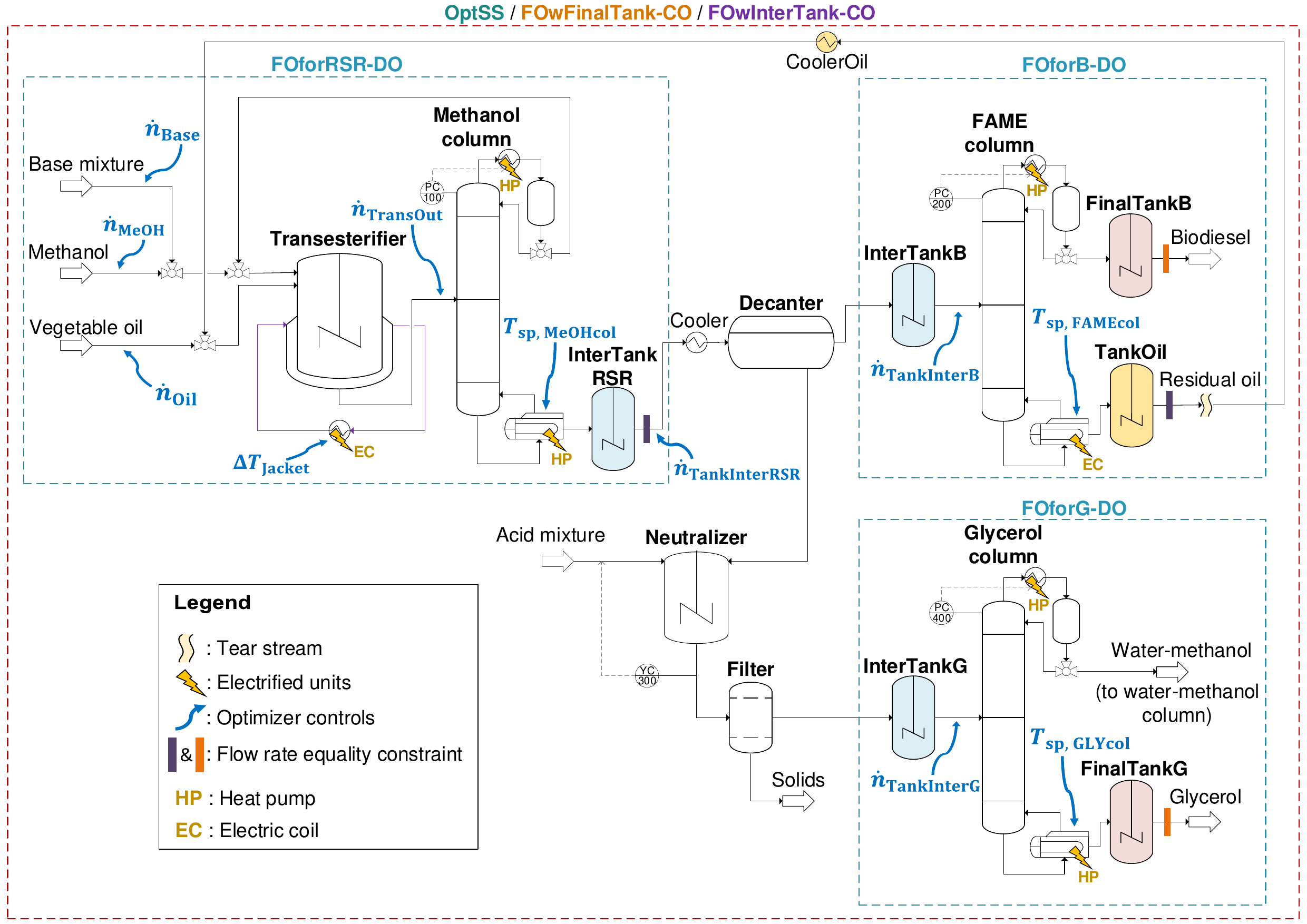}
		\caption{Superstructure flowsheet of the considered configurations of the biodiesel production process. We denote by $\dot{n}_{\text{Oil}}$, $\dot{n}_{\text{MeOH}}$, and $\dot{n}_{\text{Base}}$ the molar flow rates of the vegetable oil, methanol, and base mixture feeds, respectively. The outlet molar flow rates of the transesterifier and buffer tanks InterTankRSR, InterTankB, and InterTankG, are indicated by $\dot{n}_{\text{TransOut}}$, $\dot{n}_{\text{InterTankRSR}}$, $\dot{n}_{\text{InterTankB}}$, and $\dot{n}_{\text{InterTankG}}$, respectively. The temperature change of the transesterifier jacket medium after passing through the external heat exchanger is $\Delta T_{\text{Jacket}}$. The temperature setpoints of the column reboilers are indicated by $T_\text{sp,MeOHcol}$, $T_\text{sp,FAMEcol}$, and $T_\text{sp,GLYcol}$, respectively. All buffer tanks are excluded for the SS optimization case, \colorSS{\CentSS}, while the DO \colorCentFinal{\CentFinal} includes only the final buffer tanks FinalTankB and FinalTankG. The distributed optimization subproblems \colorDist{\DistRSR}, \colorDist{\DistB}, and \colorDist{\DistG} include all buffer tanks, whereas the DO \colorCentInter{\CentInter} exclude TankOil. The cooler, CoolerOil, is only considered for the distributed optimization case.}
		\label{fig:PFD}
	\end{figure}
	
	\subsubsection{Base-layer control system}
	
	A base-layer control system consisting of proportional-integral (PI) controllers (cf.~Figure~\ref{fig:PFD}) is implemented to determine the condenser cooling duties for the design pressure values at the top of the distillation columns by pressure controllers (PCs), and the acid feed flow rate based on a fixed pH value at the outlet of the neutralizer by a pH controller (YC). We control the liquid levels in the distillate drums, reboiler kettles, decanter, and neutralizer by cascade controllers for fixed residence times. Temperature controllers (TCs) are employed to control the reboiler temperatures to their setpoints by manipulating the reboiler heating duties. The temperature setpoints are determined by an optimizer.

	\subsubsection{Buffer tanks}
	
	We consider energy and material holdups for buffer tanks with half-full nominal liquid levels. We assume the tanks to be thermally insulated.
	
	\subsubsection{Water-methanol waste stream} \label{sec:distFit}
	
	To account for the water-methanol waste stream (cf.~Figure~\ref{fig:PFD}) in the optimization calculations, we assume the presence of a downstream purification unit that separates water and methanol (water-methanol column). The energy cost incurred by this unit, as well as the revenue generated by its purified methanol product, are considered contributors to the overall energy cost and product revenue. Therefore, based on SS simulations of a rigorous model for the water-methanol column, we fit an equation that describes the total electrical power demand of this purification unit and its methanol product flow rate as a function of the water-methanol waste flow rate and its methanol mass fraction. An equation that relates the purity of the purified methanol as a function of this mass fraction is also fitted.
 
	\subsection{Electrification of process units}
	
	The condensers and reboilers of the distillation columns, as well as the heating unit of the transesterifier, are the power-consuming components. In the following discussion, we explain how we electrify these units in this work.
	
	\paragraph{Condensers}
	We model the needed electric power for cooling in distillation column condensers through an ideal vapor-compression refrigeration cycle. We use a single-stage compression with an isentropic efficiency of $0.8$ \cite{Green.2019,Santa.2022}. We opt for ammonia as a refrigerant due to its high coefficient of performance (COP) compared to other refrigerants, and widespread use in industrial systems \cite{Baakeem.2018,McLinden.2020}. To enhance heat transfer during the refrigeration process, air fans are employed. We determine their needed electrical power by a linear correlation with the transferred heat \cite{Karakoc.2016}. The resulting COP values are around $2.7$.
	
	\paragraph{Reboilers}
	To electrify the distillation column reboilers, we utilize heat pumps or electric coils with different systems employed based on the required heat sink temperature. We model the heat pumps similarly to the refrigeration cycles used for the condensers. The methanol column reboiler operates within a temperature range of around $65$~\unit{\degreeCelsius}, for which we use a single-stage heat pump with cyclopentane as the working fluid. Such systems are characterized by high COP values for the given source and sink temperatures compared to other fluids \cite{Eppinger.2020, Frate.2019}. As a result, the COP values are around $3.5$.
	For the glycerol column reboiler, which operates within a temperature range of around $140$~\unit{\degreeCelsius}, we employ a two-stage cascade heat pump as described in \cite{Aikins.2013} and \cite{MotaBabiloni.2018}.  We use cyclopentane as the working fluid for the low-temperature stage, while water is used for the high-temperature stage, as suggested by \cite{Arpagaus.2018}. The cascade structure enables the single stages to operate at different pressure and temperature ranges, which are appropriate for the specific refrigerant. The COP values for this system are around $2$. 
	Additionally, all heat pumps are equipped with an internal heat exchange system, which reduces the mass flow of the working fluid and thereby improves the COP, as noted by \cite{Santa.2022} and \cite{MateuRoyo.2021}. Furthermore, the compression process is carried out in two steps with intercooling, which further enhances the COP and reduces the thermal stress on the compressors due to lower temperatures \cite{Aikins.2013,Cavallini.2005}.
	The temperature operating range of the FAME column reboiler is around $295$~\unit{\degreeCelsius}. We thus employ an internal electric coil with a power-to-heat efficiency of \SI{99}{\percent} as a suitable heating source for such temperature values \cite{Arpagaus.2018,Muller.2020,Schoeneberger.2022}.
	
	\paragraph{Transesterifier}
    To regulate the temperature inside the transesterifier, its jacket medium circulates through an external loop, with valves controlling the flow of hot or room-temperature water, as well as purging, depending on the desired heating or cooling mode. Heating is achieved using an electric coil, so the electrical power consumption of the transesterifier is  attributed to the heating mode only. We model this by applying complementarity constraints (CCs) as follows:
    \begin{linenomath}
	\begin{subequations}
		\begin{align}
			\dot{Q} &= \dot{Q}_{\text{h}} - \dot{Q}_{\text{c}} \,\text{,} \\
			0 \le \dot{Q}_{\text{c}} &\perp \dot{Q}_{\text{h}} \ge 0 \,\text{,} \label{eq:CC}
		\end{align}
    \end{subequations}
    \end{linenomath}
	where $\dot{Q}_{\text{c}}$ and $\dot{Q}_{\text{h}}$ are the cooling and heating power demand, respectively. The overall power demand $\dot{Q}$ depends on the process inputs and variables. We model the external loop of the transesterifier jacket medium by a heat exchanger through which the jacket medium changes its temperature by $\Delta{T_\text{Jacket}}$. We thus control the cooling or heating modes via a single control loop with one manipulated variable \cite{ElWajeh.2023}.

    Using \eqref{eq:CC} to model the transesterifier heating duty in optimization results in mathematical programs with complementarity constraints, which are challenging for nonlinear program (NLP) solvers \cite{Biegler.2010,Caspari.2020}. To allow for the use of standard NLP solvers and DAE integrators, we thus reformulate \eqref{eq:CC}, using the Fischer-Burmeister function with the smoothing term $\epsilon$ \cite{Fischer.1992,Caspari.2020}, to a smooth nonlinear complementary problem (NCP) function as follows:
    \begin{linenomath}	
    \begin{align}
		\dot{Q}_{\text{c}} + \dot{Q}_{\text{h}} \coloneqq \sqrt{\dot{Q}^2_{\text{c}} + \dot{Q}^2_{\text{h}} + \epsilon} \,\text{.} \label{eq:sNCP}
	\end{align}
    \end{linenomath}
	This NCP function is equivalent to an equality path constraint that can be incorporated directly into the integrator and solved along with the DAE system of the process model when using sequential optimization methods, unlike full discretization methods \cite{Patrascu.2018,Caspari.2020}. In this work, we use single-shooting \cite{Brusch.1973,Sargent.1978} as a direct sequential approach to solve the implemented DO problems. Without the smoothing term $\epsilon$, \eqref{eq:sNCP} leads to a nonsmooth DAE system, requiring special treatment for integration and sensitivity analysis. The value of $\epsilon$ should be sufficiently small to ensure accurate convergence of the DAE system to the exact solution, yet not overly small impeding the NLP solver's ability to explore the search space beyond the initial guess \cite{Ralph.2004,Caspari.2020}. 
	
	\section{Process optimization for flexible operation using buffer tanks} \label{sec:optimization}
	
    In this section, we present the considered process configurations and optimizations for incorporating buffer tanks to enable flexible dynamic operation, along with the benchmark standard SS operation. We first provide the general mathematical formulations of the implemented DO problems. Afterward, we discuss the approach for solving the SS optimization problem for the benchmark process design, as well as the DOs considered for different flowsheet configurations using the intermediate and final buffer tanks. For all considered optimizations, we provide the buffer tanks and process units included, process modifications, and all controls and constraints in Table~\ref{tab:optVar}.

	\subsection{Mathematical formulation}\label{subsec:mathform}
 
	Based on the modeling approach and the smooth approximation \eqref{eq:sNCP} for the transesterifier heating duty, the developed process models for all considered optimizations are smooth DAE systems of index-1. Accordingly, we solve, in all considered process configurations, DO problems on a finite time horizon $\mathcal{T} = [t_0,t_f]$ of the following form:
    \begin{linenomath}
	\begin{subequations}
		\label{eq:prob}
		\begin{align}
		\min_{\textbf{u},\textbf{x},\textbf{y}} \; \Phi(\textbf{u},\textbf{x},\textbf{y},p_{\text{el}},\textbf{v}) &= \int \limits_{t_0}^{t_f} -L(\textbf{u}(t),\textbf{x}(t),\textbf{y}(t),p_{\text{el}}(t),\textbf{v}) \mathop{\text{d}t} \,, \label{eq:probObj} \\
		s.t. \quad \textbf{M} \dot{\textbf{x}}(t) &= \textbf{f}(\textbf{u}(t),\textbf{x}(t),\textbf{y}(t),p_{\text{el}}(t),\textbf{v}), \, \forall  t \in \mathcal{T} \,, \label{eq:probf} \\
		\textbf{0} &= \textbf{g}(\textbf{u}(t),\textbf{x}(t),\textbf{y}(t),p_{\text{el}}(t),\textbf{v}), \forall t \in \mathcal{T} \,, \label{eq:probg} \\
		\textbf{0} &= \textbf{h}(\textbf{x}(t_0),\textbf{y}(t_0),p_{\text{el}}(t_0),\textbf{v}) \,, \label{eq:probh} \\
		\textbf{0} &\geq \textbf{c}(\textbf{u}(t),\textbf{x}(t),\textbf{y}(t),p_{\text{el}}(t),\textbf{v}), \, \forall t \in \mathcal{T} \,, \label{eq:probc}
		\end{align}
	\end{subequations}
    \end{linenomath}
	where $\textbf{f}: \mathcal{X} \rightarrow \mathbb{R}^{n_x}$ and $\textbf{g}: \mathcal{X} \rightarrow \mathbb{R}^{n_y}$ describe the differential-algebraic system of the process model with the non-singular and constant matrix \textbf{M} $\in \mathbb{R}^{n_x \times n_x}$, while $\mathcal{X} := \mathbb{R}^{n_x} \times \mathbb{R}^{n_y} \times \mathbb{R}^{n_u} \times \mathbb{R} \times \mathbb{R}^{n_v}$.
	The initial conditions are indicated by $\textbf{h}: \mathbb{R}^{n_x} \times \mathbb{R}^{n_y} \times \mathbb{R} \times \mathbb{R}^{n_v} \rightarrow \mathbb{R}^{n_x}$, and $\textbf{c}: \mathcal{X} \rightarrow \mathbb{R}^{n_c}$ represents the path and endpoint constraints.
	We denote the control, state, and algebraic variables by $\textbf{u}:\mathcal{T} \rightarrow \mathbb{R}^{n_u}$, $\textbf{x}:\mathcal{T} \rightarrow \mathbb{R}^{n_x}$, and $\textbf{y}:\mathcal{T} \rightarrow \mathbb{R}^{n_y}$, respectively.
	The predefined time-variant parameter, which is the electricity prices, is given by $p_{\text{el}}:\mathcal{T} \rightarrow \mathbb{R}$, while $\textbf{v}:\mathcal{T} \rightarrow \mathbb{R}^{n_v}$ are the predefined time-invariant parameters, which are the production rate demands and material prices.
	The initial and final times are denoted by $t_0 \in \mathbb{R}$ and $t_f \in \mathbb{R}$, respectively.
	The objective function $\Phi$ generally consists of the operating profit $L: \mathcal{X} \rightarrow \mathbb{R}$, but depending on the considered optimization case, $L$ can consist of the operating costs only or the power consumption instead of the cost.
	In the following sections, we define $L$ for each of the considered process configurations and the corresponding optimization problems.
	
	\subsection{Steady-state optimization via dynamic terminal-state optimization} \label{sec:optSS}
 
    We consider here the base case that we use for comparison purposes, which corresponds to a standard process design leading to a SS optimization. In this process configuration, all buffer tanks in the flowsheet in Figure~\ref{fig:PFD} are excluded. All controls, which are depicted by the blue arrows, are considered except for the InterTankRSR, InterTankB, and InterTankG outlet molar flow rates ($\dot{n}_{\text{InterTankRSR}}$, $\dot{n}_{\text{InterTankB}}$, and $\dot{n}_{\text{InterTankG}}$, respectively). Also, the CoolerOil unit and the flow rate equality constraints after InterTankRSR and TankOil are excluded here. We denote by \colorSS{\CentSS} the SS optimization case.
  
	We aim to determine an optimal constant operation benchmark that produces the required production rates of biodiesel and glycerol, and satisfies all operational requirements while minimizing power consumption and material costs and maximizing product revenues. To achieve this, we could set the right-hand side of \eqref{eq:probf} of the DAE system in \eqref{eq:prob} equal to zero and solve the resulting SS optimization problem. However, we use here an alternative approach via optimizing the dynamic terminal-state of the DAE system, which is obtained by integrating for an extended period and using constant control variables. This approach is considered more robust in the literature, e.g., \cite{Pattison.2015,Caspari.2019b}. Therefore, the controls $\textbf{u}$ in \eqref{eq:prob} are considered constant in the SS optimization. Moreover, $p_{\text{el}}$ is excluded from \eqref{eq:probObj} since the results are independent of any specific electricity price profile. In this case, $L$ is defined as:
    \begin{linenomath}
	\begin{subequations}
		\begin{align}
			L(t) &= \sum_{i=1}^{n_{\text{Prod}}} \dot{m}_i(t)v_i - qW_{\text{Tot}}(t) - \sum_{j=1}^{n_{\text{Feed}}} \dot{m}_j(t)v_j \, \text{,} \label{eq:L_SS} \\
			W_{\text{Tot}}(t) &= W_{\text{Trans}}(t) + W_{\text{Mc}}(t) + W_{\text{Fc}}(t) + W_{\text{Gc}}(t) + W_{\text{WMc}}(t) \, \text{,}
		\end{align}
	\end{subequations}
    \end{linenomath}
	where $\dot{m}_i$ and $\dot{m}_j$ indicate the production and consumption rates of products and feeds, with the corresponding specific material prices, denoted by $v_i$ and $v_j$, respectively. The $n_{\text{Prod}}$ products are biodiesel, glycerol, solids, and purified methanol by the water-methanol column (cf.~Figure~\ref{fig:PFD} and Section~\ref{sec:distFit}). The $n_{\text{Feed}}$ feeds are vegetable oil, methanol, base mixture, and acid mixture. The total power demand is given by $W_{\text{Tot}}$, and weighted by the factor $q$. The power demands of the transesterifier, methanol, FAME, glycerol, and water-methanol columns are indicated by $W_{\text{Trans}}$, $W_{\text{Mc}}$, $W_{\text{Fc}}$, $W_{\text{Gc}}$, $W_{\text{WMc}}$, respectively.
	
	We use a final time $t_f$ of two days, which is sufficiently large for the DAE system to obtain a new quasi SS, starting at the initial DAE state defined by \eqref{eq:probh}. We selected the time horizon by forward simulation of the model, ensuring that a SS is obtained. By optimizing the terminal-state of the DAE system using constant controls, we obtain one feasible SS solution. However, there may be multiple solutions, and the found one may not necessarily be stable. Therefore, we perform a stability check by linearizing the DAE system at the found SS solution and using the indirect method of Lyapunov. We observe that all the real parts of the eigenvalues of the linearized system matrix are negative, indicating that the DAE system is stable at the found SS solution.
	
	\subsection{Process configurations and dynamic optimization for flexible operation}
	
	To enable flexible process operation, we add final and/or intermediate buffer tanks. We aim to dynamically operate the process by solving \eqref{eq:prob} for the different buffer tank configurations while producing the same amount of biodiesel and glycerol as the SS operation benchmark within the considered time horizon. The controls $\textbf{u}$ are time-variant variables, and $p_{\text{el}}$ is kept in \eqref{eq:probObj}. Accordingly, the operating profit $L$ is defined as follows:
    \begin{linenomath}
	\begin{align}
		L(t)=  \sum_{i=1}^{n_{\text{Prod}}} \dot{m}_i(t)v_i - p_{\text{el}}(t) W_{\text{Tot}}(t) - \sum_{j=1}^{n_{\text{Feed}}} \dot{m}_j(t)v_j \, \text{.} \label{eq:L_Dyn}
	\end{align}
    \end{linenomath}
	First, we discuss the configuration of adding final buffer tanks exclusively, before moving on to the use of intermediate tanks. Additionally, flexibility potential exists not only in production rates but also in the purity of final products. In Section~\ref{sec:flexPurity}, we elaborate on how product purity specification can be used to achieve more process flexibility, particularly for the glycerol product.
	
	\subsubsection{Process flexibilization via final buffer tanks only}
	
	As shown in Figure~\ref{fig:PFD} and Table~\ref{tab:optVar}, we include in this case only the two final buffer tanks, namely FinalTankB and FinalTankG, in the process flowsheet. All intermediate buffer tanks are excluded here. Final buffer tanks are essential when aiming for flexible operation while simultaneously meeting specific production demands, such as biodiesel and glycerol production in our study. We solve \eqref{eq:prob} for this configuration, investigating the operational flexibility that can be achieved through the use of final buffer tanks only. The operating profit $L$ is defined here as in \eqref{eq:L_Dyn}. We represent the optimization for this case as \colorCentFinal{\CentFinal} (flexible operation with final tanks - centralized optimization).
	
	We impose additional endpoint constraints for the liquid levels in FinalTankB and FinalTankG (cf.~Table~\ref{tab:optVar}) to guarantee that the optimizer does exploit the initial holdups in the tanks to satisfy production demands. Therefore, the liquid levels have to be maintained at their initial values at the end of the time horizon. Our analysis aims to investigate the production flexibility of both biodiesel and glycerol products while ensuring that their respective demands are met.
	
	\subsubsection{Process flexibilization via final and intermediate buffer tanks}
	
	Various unit operations have different potentials for flexibility, depending on their operational requirements and positions in the process. For instance, a unit operation's ability to be flexible in production, and thus in its power consumption, may be limited by downstream processes. In particular, liquid level limits in downstream processes may impede the production flexibility potentials of upstream processes. As a solution, the incorporation of additional intermediate buffer tanks can decouple the dynamics between process parts and consequently render full exploitation of the production flexibility of the overall process.
	
	In this work, the output production rate of the methanol column bottom is restricted by the liquid level limits in the downstream processes, notably, the decanter and columns. Thus, we incorporate the buffer tank InterTankRSR to realize the full potential of production flexibility for the methanol column (cf.~Figure~\ref{fig:PFD} and Table~\ref{tab:optVar}). Consequently, the outlet flow rate of InterTankRSR is a new control variable in the DO problem. 
	
	Furthermore, to enable varying the production rates through the other downstream power-consuming units, namely the FAME and glycerol columns, we include the buffer tanks, InterTankB and InterTankG, respectively. As a result, the outlet flow rates of these tanks need to be manipulated, making them new controls in the DO problem (cf.~Figure~\ref{fig:PFD} and Table~\ref{tab:optVar}). For all intermediate buffer tanks, we add endpoint constraints for their liquid levels as well as for the purities of their content species (cf.~Table~\ref{tab:optVar}) to ensure that the state at the final time is equal to its initial value.
	
	In addition to fully exploiting production flexibility, the added intermediate buffer tanks can also facilitate the implementation of multiple distributed optimizers for different process parts, rather than relying on a single centralized optimizer for the entire process. Given the high computational cost and convergence challenges associated with large-scale DO problems, we propose two optimization strategies for employing intermediate buffer tanks. First, we examine the application of a centralized monolithic optimizer for the entire process. Subsequently, we introduce an additional buffer tank for the residual oil recycle (cf.~Figure~\ref{fig:PFD} and Table~\ref{tab:optVar}), and employ three distributed optimizers instead.
	
	\paragraph{Centralized monolithic optimization} We denote the optimizer that solves \eqref{eq:prob}, where the operating profit $L$ is defined as in \eqref{eq:L_Dyn}, for the entire process including all intermediate buffer tanks by \colorCentInter{\CentInter} (flexible operation with intermediate and final tanks - centralized optimization). In this case, the convergence of the NLP solver is sensitive to the initial guess and the scaling of controls and constraints. In DO problems with large-scale DAE systems, the optimization algorithm is particularly susceptible to encountering ill-conditioning issues, which adversely impact convergence. Utilizing distributed optimizers for different parts of the process leads to DO problems with smaller DAE systems, thus, reducing the number of variables that the user needs to initialize and scale. Consequently, the optimization algorithm is less prone to non-convergence and is computationally more efficient.
	
	\paragraph{Distributed optimization} By introducing an additional buffer tank, TankOil, and the water-operating cooler, CoolerOil, for the residual oil recycle, we fix its flow rate and temperature, thereby, allowing to tear this recycle stream (cf.~Figure~\ref{fig:PFD} and Table~\ref{tab:optVar}). Consequently, we can decouple the upstream processes of the InterTankRSR tank from its downstream processes. However, to enable this, we need additional constraints. Specifically, we constrain the outlet flow rates of InterTankRSR and TankOil tanks to fixed flow rates and impose path constraints on the species purities of the InterTankRSR outlet (cf.~Table~\ref{tab:optVar}). Therefore, we can implement the three optimizers \colorDist{\DistRSR}, \colorDist{\DistB}, and \colorDist{\DistG}, as shown in Figure~\ref{fig:PFD}. Collectively, we refer to them as \colorDistAll{\DistAll} (flexible operation with intermediate and final tanks - distributed optimization). This approach involves solving three DO problems of \eqref{eq:prob}, with smaller DAE systems, leading to fewer non-convergence issues and less computational cost. However, it is important to note that these additional constraints are restrictive. Notably, fixing the flow rate of the residual oil recycle stream results in less efficient material consumption of the oil feed. This increases the overall material costs compared to the \colorCentInter{\CentInter} case where we have full degrees of freedom. In Section~\ref{sec:results}, we examine whether these restrictions are significant. The operating profits for the three problems are given by:
    \begin{linenomath}
	\begin{subequations}
	\begin{align}
		L_{\text{\DistRSR}}(t) &= - p_{\text{el}}(t) \left(W_{\text{Trans}}(t) + W_{\text{Mc}}(t)\right) - \sum_{j=1}^{n_{\text{Feed, \DistRSR}}} \dot{m}_j(t)v_j \, \text{,} \label{eq:L_\DistRSR} \\
		L_{\text{\DistB}}(t) &= \dot{m}_{\text{B}}(t)v_{\text{B}} - p_{\text{el}}(t) W_{\text{Fc}}(t) \, \text{,} \label{eq:L_\DistB} \\
		L_{\text{\DistG}}(t) &= \sum_{i=1}^{n_{\text{Prod, \DistG}}} \dot{m}_i(t)v_i - p_{\text{el}}(t) \left(W_{\text{Gc}}(t) + W_{\text{WMc}}(t)\right) \, \text{,} \label{eq:L_\DistG}
	\end{align}
	\end{subequations}
    \end{linenomath}
	such that the acid feed is not included in $L_{\text{\DistRSR}}$; ${\text{B}}$ indicates biodiesel; and the products in $L_{\text{\DistG}}$ are glycerol, solids, and purified methanol by the water-methanol column.

    It is worth noting that the additional TankOil serves the sole purpose of facilitating the decoupling of process parts and enabling distributed optimizations. Consequently, it does not provide any additional flexibility benefits for the centralized monolithic optimization case. In our process, the two recycle streams from the methanol and FAME columns are already buffered within the transesterifier, rendering intermediate buffer tanks unnecessary for these streams.

	\subsubsection{Flexible purity production} \label{sec:flexPurity}
	
	Enforcing the purity constraints of the final products at the outlet streams of the buffer tanks, FinalTankB and FinalTankG, enables the flexibilization of the produced product purities at their inlet streams. Specifically, a higher degree of freedom is given to the optimizer to vary the purity at the buffer tank inlet streams while satisfying its requirements at the outlet side. Producing higher purities in distillation columns is associated with higher power consumption. As a result, flexibilizing the produced purities of biodiesel and glycerol products leads to additional savings in energy costs. Thus, we impose the required purity limits for both biodiesel and glycerol products at the outlet streams of the corresponding final buffer tanks. However, flexibility in purity production is more significant when having purity limits for a single species only. That is the case for the glycerol product in our study, where we only control glycerol species purity. In contrast, there are, practically, purity limits for FAME, methanol, water, glycerol, and monoolein for the biodiesel product.

	\begin{table}[htbp]
		\centering
		\caption{Summary of the included buffer tanks and units, process modifications, control variables, and operational constraints for the implemented process configurations. Due to the thermal degradation limits of biodiesel and glycerol products, the maximum temperatures in the FAME and glycerol column reboilers are $300$~\unit{\degreeCelsius} and $150$~\unit{\degreeCelsius}, respectively. The maximum temperature changes of the transesterifier jacket medium $\Delta T_\text{Jacket}$ are limited to $\pm10$~\unit{\degreeCelsius}. Time-variant is indicated by TV. We denote by LL and EC the liquid levels and equality constraints, respectively. Purities are in \unit{\kilogram\per\kilogram}.}
        {\fontsize{11pt}{11pt}\selectfont
		\begin{tabular}[t]{>{\raggedright}m{0.52\linewidth}>{\centering}m{0.06\linewidth}>{\centering}m{0.09\linewidth}>{\centering}m{0.13\linewidth}>{\centering\arraybackslash}m{0.1\linewidth}}
			\toprule
			& \colorSS{\CentSS} & \colorCentFinal{FOwFinal-Tank-CO} & \colorCentInter{FOwInter-Tank-CO} & \colorDistAll{FOwInter-Tank-DO} \\
			\midrule
			Additional process units and modifications & & & & \\
			\midrule
			FinalTankB and FinalTankG & --- & \checkmark & \checkmark & \checkmark \\
			InterTankRSR, InterTankB, and InterTankG & --- & --- & \checkmark & \checkmark \\
			TankOil and CoolerOil & --- & --- & --- & \checkmark \\
			Tearing the residual oil recycle stream & --- & --- & --- & \checkmark \\
			\midrule
			Controls & Constant & TV & TV & TV \\
			\midrule
			$\dot{n}_\text{Oil} \; \left[\unit{\kilo\mol\per\hour} \right]$ & \checkmark & \checkmark & \checkmark & \checkmark \\
			$\dot{n}_\text{MeOH} \; \left[\unit{\kilo\mol\per\hour} \right]$ & \checkmark & \checkmark & \checkmark & \checkmark \\
			$\dot{n}_\text{Base} \; \left[\unit{\kilo\mol\per\hour} \right]$ & \checkmark & \checkmark & \checkmark & \checkmark \\
			$\dot{n}_\text{TransOut} \; \left[\unit{\kilo\mol\per\hour} \right]$ & \checkmark & \checkmark & \checkmark & \checkmark \\
			$\Delta T_\text{Jacket} \; \left[\unit{\kelvin} \right]$, limited to $\pm10$~\unit{\degreeCelsius}& \checkmark & \checkmark & \checkmark & \checkmark \\
			$T_\text{sp,MeOHcol} \; \left[\unit{\kelvin} \right]$ & \checkmark & \checkmark & \checkmark & \checkmark \\
			$T_\text{sp,FAMEcol} \; \left[\unit{\kelvin} \right]$, upper bound of $300$~\unit{\degreeCelsius} & \checkmark & \checkmark & \checkmark & \checkmark \\
			$T_\text{sp,GLYcol} \; \left[\unit{\kelvin} \right]$, upper bound of $150$~\unit{\degreeCelsius} & \checkmark & \checkmark & \checkmark & \checkmark \\
			$\dot{n}_\text{InterTankRSR} \; \left[\unit{\kilo\mol\per\hour} \right]$ & --- & --- & \checkmark & --- \\
			$\dot{n}_\text{InterTankB} \; \left[\unit{\kilo\mol\per\hour} \right]$ & --- & --- & \checkmark & \checkmark \\
			$\dot{n}_\text{InterTankG} \; \left[\unit{\kilo\mol\per\hour} \right]$ & --- & --- & \checkmark & \checkmark \\
			\midrule
			Path and endpoint constraints & & & & \\
			\midrule
			Transesterifier LL $\left[\unit{\meter} \right]$ & \checkmark & \checkmark & \checkmark & \checkmark \\
			FinalTankB and FinalTankG LL $\left[\unit{\meter} \right]$ & --- & \checkmark & \checkmark & \checkmark \\
			InterTankRSR, InterTankB, and InterTankG LL~$\left[\unit{\meter} \right]$ & --- & --- & \checkmark & \checkmark \\
			TankOil LL $\left[\unit{\meter} \right]$ & --- & --- & --- & \checkmark \\
			InterTankRSR content purities & --- & --- & --- & \checkmark \\
			\midrule
			Path constraints & & & & \\
			\midrule
			Decanter and neutralizer LL $\left[\unit{\meter} \right]$ & \checkmark & \checkmark & \checkmark & --- \\
            Columns:~reboilers,~distillate~drums,~and~trays~LL~$\left[\unit{\meter} \right]$ & \checkmark & \checkmark & \checkmark & \checkmark \\
			EN 14214 \cite{BritishStandardsInstitution.2010} biodiesel purities & \checkmark & \checkmark & \checkmark & \checkmark \\
			$99$~wt\unit{\percent} glycerol purity & \checkmark & \checkmark & \checkmark & \checkmark \\
			Biodiesel production demand $\left[\unit{\kilogram\per\hour} \right]$ (EC) & \checkmark & \checkmark & \checkmark & \checkmark \\
			Glycerol production demand $\left[\unit{\kilogram\per\hour} \right]$ (EC) & \checkmark & \checkmark & \checkmark & \checkmark \\
			InterTankRSR outlet $\left[\unit{\kilogram\per\hour} \right]$ (EC) & --- & --- & --- & \checkmark \\
			TankOil outlet $\left[\unit{\kilogram\per\hour} \right]$ (EC) & --- & --- & --- & \checkmark \\
			\midrule
			Endpoint constraints & & & & \\
			\midrule		
			Transesterifier content purities & --- & \checkmark & \checkmark & \checkmark \\
			InterTankRSR, InterTankB, and InterTankG content purities & --- & --- & \checkmark & \checkmark \\
			\bottomrule
		\end{tabular}
        }
		\label{tab:optVar}
	\end{table}

	\section{Operational scenario} \label{sec:scenario}

\begin{figure}[htbp]
		\centering
		\includegraphics[width=0.25\linewidth]{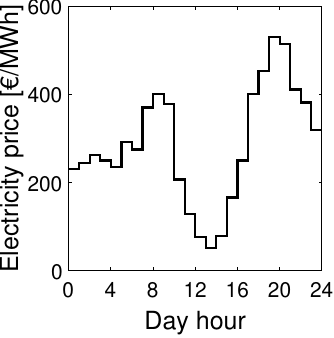}
		\caption{German day-ahead electricity prices for September~3,~2022 \cite{SMARD}.}
		\label{fig:Scenario}
	\end{figure}

	Here, we introduce the considered demand-response scenario for all optimization cases, that is, how we define $p_{el}(t)$ and $\textbf{v}$. We conduct simulations over a time horizon of one day, during which constant production demands of $20$~\unit{\tonne\per\hour} for biodiesel and $2.12$~\unit{\tonne\per\hour} for glycerol are required. We use historical electricity price data from the German day-ahead spot market for September~3,~2022 \cite{SMARD}, which is depicted in Figure~\ref{fig:Scenario}. We use constant prices for raw materials and final products. For all DO problems, we discretize the control variables at an equidistant interval of one hour, while the constraints are discretized at 30-minute intervals.

	\section{Implementation} \label{sec:implementation}

	We solve all optimization problems to local NLP convergence with direct single-shooting \cite{Brusch.1973,Sargent.1978} using our open-source framework DyOS (Dynamic Optimization Software) \cite{Caspari.2019}. Using NIXE (NIXE Is eXtrapolated Euler) \cite{Hannemann.2010} as a DAE integrator and SNOPT (Sparse Nonlinear OPTimizer) \cite{Gill.2005} as an NLP, the DO problems are solved sequentially. The Modelica model is coupled to DyOS as a Functional Mockup Unit (FMU) \cite{FMI} generated by Dymola \cite{DassaultSystemes.2023}. As FMU only supports ODEs, Dymola performs numerical reduction and symbolic reformulation of the DAE system to provide an FMU. In addition, we use Dymola for model linearization and for calculating the eigenvalues of the linearized system matrix. We set the DAE integrator, NLP feasibility, and optimality tolerances~to~$10^{-4}$.
	
	\section{Results and discussion} \label{sec:results}
	
	We present production rates and the total power demand results for all the considered process configurations first, before we discuss how the buffer tanks are utilized to enhance the production and, thus, power consumption flexibility. Afterward, we demonstrate how glycerol product purity can vary based on energy costs. We also evaluate the economic performance of the considered optimizations and compare them with the SS operation benchmark \colorSS{\CentSS}. Lastly, we compare the computational performance of \colorDistAll{\DistAll} and \colorCentInter{\CentInter}.
	
	\subsection{Production rate and power demand} \label{sec:prodPow}
	
\begin{figure}[htbp]
		\centering
		\includegraphics[width=0.95\linewidth]{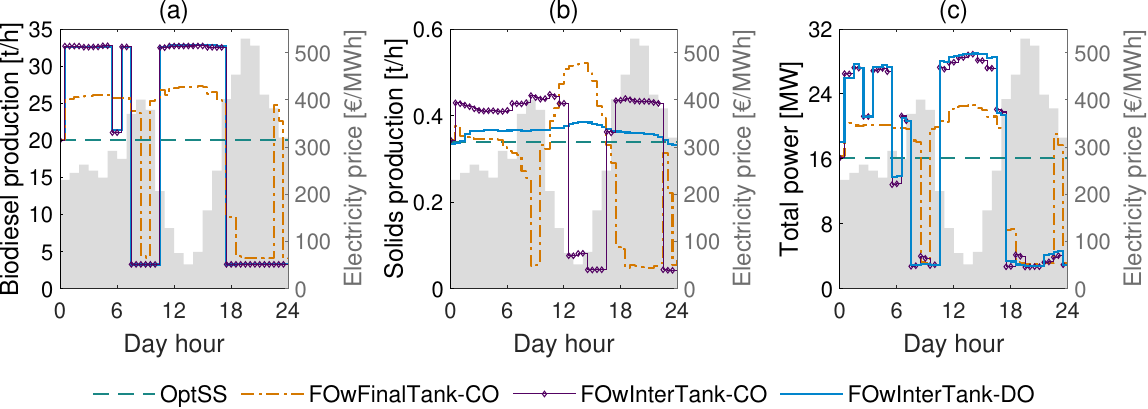}
		\caption{Production rates and total power demand for all optimizers. (a): biodiesel production rates; (b): solids production rates; (c): total power consumption rates. The shaded areas illustrate the electricity price profile, which corresponds to the secondary axes.}
		\label{fig:ProdPow}
	\end{figure}

	The biodiesel production rates are presented in Figure~\ref{fig:ProdPow}a for all considered process configurations. We observe that all DO optimizations, as compared to \colorSS{\CentSS}, enable production flexibility based on electricity prices. We see how the production rate profiles are opposite to that of the electricity prices. The DO optimizers promptly increase the production rates to move away from the nominal starting point due to lower electricity prices during this period. Subsequently, the production rates are adjusted according to prices to optimize the operating profit while adhering to operational constraints, particularly level limits. The intermediate buffer tanks provide significant additional flexibility, as evident from the profiles of \colorCentInter{\CentInter} and \colorDistAll{\DistAll} compared to that of \colorCentFinal{\CentFinal}. Notably, in our case, the upper limits of production rates are increased, e.g., during the period between $10$ and $17$. On the other hand, \colorCentFinal{\CentFinal} cannot fully leverage the low prices or always maintain production rates at minimum levels during high prices. During periods around $8$, $10$, and $23$, \colorCentFinal{\CentFinal} operates at high production rates, even exceeding the nominal rate, despite high prices. In contrast, \colorCentInter{\CentInter} and \colorDistAll{\DistAll} maintain the production rate at its minimum operating limit. In addition, the production in \colorDistAll{\DistAll} is identical to that of \colorCentInter{\CentInter}, indicating that the use of distributed optimizers rather than a centralized monolithic has a similar flexibility potential. Moreover, by comparing the profiles of \colorCentFinal{\CentFinal} with those of \colorCentInter{\CentInter} and  \colorDistAll{\DistAll}, it is evident that the dynamics in the latter cases are buffered due to the intermediate buffer tanks.
	
	The production rate profiles for the glycerol product and water-methanol waste are similar to that of biodiesel. However, during periods of increased production (e.g., between $10$ and $17$), more water enters the glycerol purification section for \colorDistAll{\DistAll} compared to the other DO cases. This is due to the additional purity constraints on the InterTankRSR outlet stream that are needed when separating the upstream processes of InterTankRSR from the downstream ones. Therefore, during periods of increased production, higher amounts of water-methanol waste are produced for \colorDistAll{\DistAll} compared to \colorCentFinal{\CentFinal} and \colorCentInter{\CentInter}, resulting in higher power demands for both glycerol and water-methanol columns.
	
	Figure~\ref{fig:ProdPow}b shows the production rates of solids, the product obtained from the filter unit. Although the produced amount is relatively insignificant compared to the primary products, it provides insights into how the buffer tanks affect the process dynamics, especially the InterTankRSR. For \colorCentFinal{\CentFinal}, the production rate is similar to that of biodiesel since this case only considers the final buffer tanks. However, for \colorCentInter{\CentInter}, the production rate is the opposite and follows the electricity price profile. The solids product is a downstream product directly after the InterTankRSR, without any buffer tanks in between. In \colorCentInter{\CentInter}, all intermediate buffer tanks have full degrees of freedom, enabling InterTankRSR to deliver high flow rates to the downstream buffer tanks, InterTankB, and InterTankG, during periods of high prices, and vice versa. On the other hand, InterTankB and InterTankG operate oppositely as they are responsible for achieving flexible production rates for biodiesel and glycerol, respectively, as shown in Figure~\ref{fig:ProdPow}a. Therefore, before the prices increase, InterTankRSR fills these tanks, explaining the high flow rates of the solids product during high prices. During low-price periods, the outlet stream of InterTankRSR decreases, allowing InterTankRSR to refill while InterTankB and InterTankG are not utilized as much as during high prices. Section~\ref{sec:tanks} provides further elaboration on the buffer tank utilization. Conversely, for \colorDistAll{\DistAll}, the solids production rate is mostly unaffected by the flexibilization of the other process parts. When using distributed optimizers and fixing the outlet flow rate of InterTankRSR, all process units between InterTankRSR and its downstream buffer tanks, InterTankB and InterTankG, operate almost constantly. Slight variations in the production rates occur due to changes in the species purities of the InterTankRSR outlet. As those species purities are not fixed, the outlet streams from the decanter vary according to their degree of separation. It is noteworthy how distributed optimizers can decouple the dynamics of various process parts using intermediate buffer tanks, allowing different parts of the process to operate at varying degrees of flexibility.
	
	The total power consumption rates are illustrated in Figure~\ref{fig:ProdPow}c, which resemble the biodiesel production rate profiles. Additionally, we see how \colorCentInter{\CentInter} and \colorDistAll{\DistAll} decrease the power consumption rates at periods around $3$ and $6$, as compared to \colorCentFinal{\CentFinal}. This again demonstrates the effectiveness of intermediate buffer tanks in exploiting slight variations in prices.
		
	\subsection{Buffer tank levels} \label{sec:tanks}
	
\begin{figure}[htbp]
		\centering
		\includegraphics[width=0.95\linewidth]{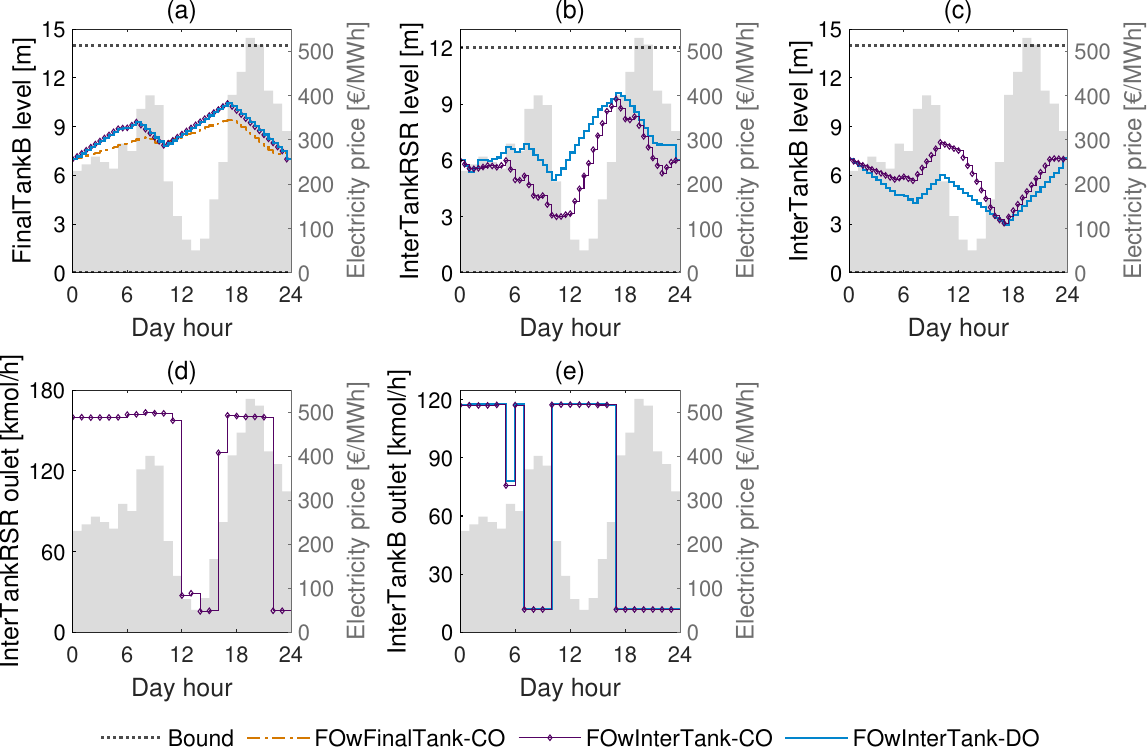}
		\caption{(a), (b), and (c) show the liquid levels in FinalTankB, InterTankRSR, and InterTankB, respectively. The control variables for InterTankRSR and InterTankB, which are the outlet flow rates, are shown in (d) and (e), respectively. For \colorDistAll{\DistAll}, the InterTankRSR outlet flow rate is an equality constraint. Therefore, only \colorCentInter{\CentInter} has it as a control variable. The secondary axes correspond to the electricity prices, which are depicted in shaded areas.}
		\label{fig:Tanks}
	\end{figure}

	The liquid levels and controls of InterTankRSR, InterTankB, and FinalTankB in Figure~\ref{fig:Tanks} demonstrate how the final and intermediate buffer tanks enable high degrees of production flexibility. In Figure~\ref{fig:Tanks}a, FinalTankB levels increase when electricity prices are low and decrease otherwise. The tank holdup is utilized during high-price periods to ensure that the required biodiesel demand is met, despite the decrease in production rates (Figure~\ref{fig:ProdPow}a). The changes in the levels are steeper for \colorCentInter{\CentInter} and \colorDistAll{\DistAll} than for \colorCentFinal{\CentFinal}, highlighting the additional flexibility that intermediate buffer tanks provide.

	For the glycerol product, the profiles for FinalTankG resemble those of FinalTankB. In contrast, InterTankB or InterTankG exhibit the opposite behavior of final tanks. These intermediate tanks supply the required flow rates for downstream units to operate at the desired capacity. Therefore, they are utilized during low-price periods and filled during high-price periods. For example, during the period between $10$ and $17$, InterTankB is used to operate at high outlet flow rates (Figure~\ref{fig:Tanks}e) so that the FAME column operates at maximum capacity. Consequently, the levels in InterTankB decrease while they increase in FinalTankB.
	
	Unlike InterTankB or InterTankG, InterTankRSR is utilized similarly to the final tanks. It is filled during low-price periods and utilized otherwise (Figure~\ref{fig:Tanks}b). InterTankRSR fills InterTankB and InterTankG before low-price periods so that the latter tanks are utilized during high-price periods. This is evident from the InterTankRSR control variable profile (only in the \colorCentInter{\CentInter} case), which shows high flow rate values during high-price periods and vice versa (Figure~\ref{fig:Tanks}d).
	
	\subsection{Flexible purity production}
			
\begin{figure}[htbp]
		\centering
		\includegraphics[width=0.95\linewidth]{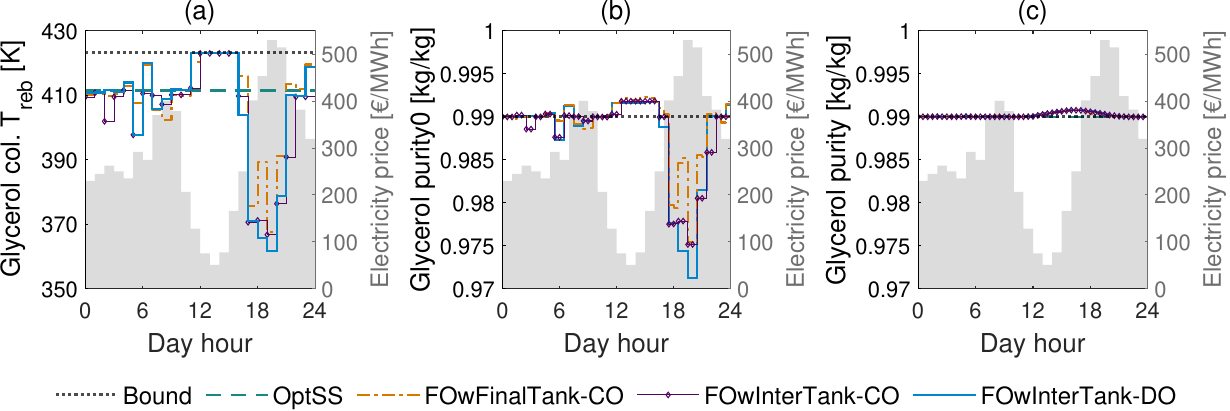}
		\caption{(a) shows the temperature setpoint of the glycerol column reboiler, which is an optimization control variable. (b) and (c) illustrate the glycerol purity in the inlet and outlet streams of FinalTankG, respectively. It can be observed from (b) that lower purities are produced during high-price periods and vice versa, while the delivered product purity (in~(c)) meets the minimum limit of $99$~wt\unit{\percent}. The profile of the electricity prices is represented by the shaded area, which corresponds to the secondary axis.}
		\label{fig:GlyPurity}
	\end{figure}

	Figure~\ref{fig:GlyPurity} shows the variability of the produced purity of the glycerol product, which enables additional operational flexibility. The figure depicts the purities before and after FinalTankG, as well as the reboiler temperature setpoint of the glycerol column, which is an optimization control variable. The temperature setpoints increase during low-price periods and decrease otherwise for all DO cases (Figure~\ref{fig:GlyPurity}a), resulting in added power demand flexibility. During periods of lowest price, specifically between $10$ and $17$, the temperature setpoints reach their maximum allowable limit (glycerol thermal degradation limit) to maximize the power demands of the glycerol column during those periods. This leads to higher glycerol purities being produced and stored in FinalTankG (Figure~\ref{fig:GlyPurity}b). On the other hand, lower purities are produced during high-price periods and mixed with the higher-purity content in FinalTankG, resulting in glycerol being delivered at the tank outlet with purities above the required limit (Figure~\ref{fig:GlyPurity}c). The purity profiles at the inlet of FinalTankG (Figure~\ref{fig:GlyPurity}b) mirror the reboiler temperature profiles, particularly how they decrease and then increase during the period of highest prices, i.e., between $17$ and $24$. During this period, the purities decrease to below $97.5$~wt\unit{\percent} (for \colorDistAll{\DistAll}), while they remain above the minimum limit at the FinalTankG outlet (Figure~\ref{fig:GlyPurity}c).

	Furthermore, we observe that the purities at the FinalTankG outlet around $17$ are not at the minimum limit (Figure~\ref{fig:GlyPurity}c), in contrast to other periods. Just before this period, i.e., between $10$ and $17$, the purities at the inlet reach their highest operational limit (due to temperature limit), resulting in slightly higher purities at the outlet during the period around $17$. Consequently, during the next period when the prices are high (between $17$ and $24$), the tank inlet purities can be minimized while maintaining the outlet purities above the minimum limit.
	
	\subsection{Economic evaluation}
	
	In Table~\ref{tab:econ}, we present the total operating profits, energy costs, and material costs for all optimization case studies. The energy costs for three parts of the process are provided: RSRprocess, which includes the transesterifier and methanol column, Bprocess for the FAME column, and Gprocess for the glycerol and water-methanol columns. Additionally, all costs and profits for the DOs are given relative to the SS optimization benchmark.

    \begin{table}[htbp]
		\small
		\centering
		\caption{Total operating profit, and energy and material costs for each optimizer. Energy costs for different parts of the process are given. RSRprocess energy includes the transesterifier and the methanol column power demands. Bprocess energy indicates the FAME column power demand, while Gprocess energy indicates that of the glycerol and downstream water-methanol columns. The economic performance (savings) relative to the SS operation case, \colorSS{\CentSS}, is given in parentheses.}
        {\fontsize{11pt}{11pt}\selectfont
		\begin{tabular}[t]{>{\raggedright}m{0.29\linewidth}>{\centering}m{0.09\linewidth}>{\centering}m{0.16\linewidth}>{\centering}m{0.16\linewidth}>{\centering\arraybackslash}m{0.16\linewidth}}
			\toprule
			& \colorSS{\CentSS} & \colorCentFinal{\CentFinal} & \colorCentInter{\CentInter} & \colorDistAll{\DistAll} \\
			\midrule
			RSRprocess energy cost [\unit{\kilo\euro}] & $11.4$ & $10.4$ (\SI{9.3}{\percent}) & $9$ (\SI{22}{\percent}) & $8.2$ (\SI{28}{\percent}) \\
			Bprocess energy cost [\unit{\kilo\euro}] & $77.8$ & $64.6$ (\SI{17}{\percent}) & $56.4$ (\SI{28}{\percent}) & $56.5$ (\SI{27}{\percent}) \\
			Gprocess energy cost [\unit{\kilo\euro}] & $22.3$ & $14.4$ (\SI{36}{\percent}) & $14.4$ (\SI{36}{\percent}) & $16.1$ (\SI{28}{\percent}) \\
			Total energy cost [\unit{\kilo\euro}] & $111.5$ & $89.4$ (\SI{20}{\percent}) & $79.8$ (\SI{29}{\percent}) & $80.8$ (\SI{28}{\percent}) \\
			Total material cost [\unit{\kilo\euro}] & $563.9$ & $556.8$ (\SI{1.3}{\percent}) & $556.3$ (\SI{1.4}{\percent}) & $562.9$ (\SI{0.2}{\percent}) \\
			Total profit [\unit{\kilo\euro}] & $639.8$ & $667$ (\SI{4.3}{\percent}) & $678.6$ (\SI{6.1}{\percent}) & $672.4$ (\SI{5.1}{\percent}) \\
			\bottomrule
		\end{tabular}
        }
		\label{tab:econ}
	\end{table}
	
	We find that \colorCentFinal{\CentFinal} incurs \SI{20}{\percent} less total energy cost relative to \colorSS{\CentSS}, while \colorCentInter{\CentInter} and \colorDistAll{\DistAll} result in \SI{29}{\percent} and \SI{28}{\percent} less total energy cost, respectively. The total energy cost savings are similar for the implementation of distributed optimizers \colorDistAll{\DistAll} compared to the centralized monolithic \colorCentInter{\CentInter}. However, there are differences in the savings for different process parts, particularly for the RSRprocess and Gprocess. We observe that tearing the recycle oil stream and fixing its flow rate, while imposing purity limits on the InterTankRSR outlet and also fixing its flow rate, leads to less power consumption in the methanol column for \colorDistAll{\DistAll} compared to \colorCentInter{\CentInter}. However, as explained in Section~\ref{sec:prodPow}, more water entering the glycerol purification process for \colorDistAll{\DistAll} results in additional power consumption in the glycerol and water-methanol columns. This explains the higher energy cost incurred for \colorDistAll{\DistAll} compared to \colorCentInter{\CentInter} and \colorCentFinal{\CentFinal}.
	
	A comparison of the RSRprocess energy costs for \colorCentFinal{\CentFinal} with that of \colorDistAll{\DistAll} and \colorCentInter{\CentInter} reveals significant differences in terms of savings, emphasizing the additional operational flexibility that InterTankRSR offers for the methanol column. In the case of Bprocess, the savings for both \colorCentInter{\CentInter} and \colorDistAll{\DistAll} are significantly higher than \colorCentFinal{\CentFinal}, indicating the added flexibility of InterTankRSR when combined with InterTankB. It is worth noting that the intermediate tanks do not lead to additional energy savings for the Gprocess, as demonstrated by the comparison of energy costs between \colorCentFinal{\CentFinal} and \colorCentInter{\CentInter}.
	
	We find only slight savings in material costs for the DO cases compared to the SS benchmark, which is expected since the material prices remain unchanged. However, when comparing \colorDistAll{\DistAll} to \colorCentFinal{\CentFinal} or \colorCentInter{\CentInter}, we observe fewer material cost savings for \colorDistAll{\DistAll} due to the restrictive consideration of tearing the recycle oil stream and fixing its flow rate. Fixing the recycle stream results in less efficient utilization of the expensive oil raw material, which is also the process throughput manipulator. Thus, the additional material cost incurred for \colorDistAll{\DistAll} can be explained by the less efficient use of the oil feed.
	
	After analyzing the total operating profits, it is evident that material costs exceed energy costs, resulting in considerably lower additional profits by an order of magnitude, as compared to energy costs savings. We also observe that \colorDistAll{\DistAll} results in \SI{1}{\percent} less additional profit (relative to \colorSS{\CentSS}) compared to \colorCentInter{\CentInter}. In Section~\ref{sec:solTime}, we compare the computational costs and comment on the NLP solver convergence for \colorCentInter{\CentInter} and \colorDistAll{\DistAll} assessing the significance of this \SI{1}{\percent} reduction in savings.
	
	\subsection{Solution times} \label{sec:solTime}
	
\begin{figure}[htbp]
		\centering
		\includegraphics[width=1\linewidth]{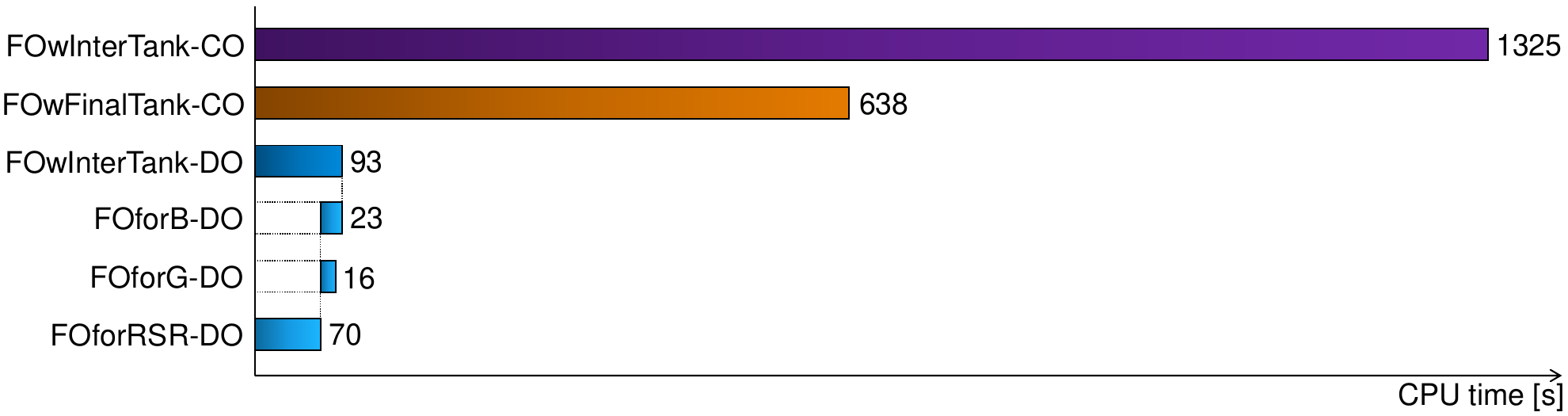}
		\caption{CPU times for solving the DO problems of the considered optimizers. After solving \colorDist{\DistRSR}, \colorDist{\DistB} and \colorDist{\DistG} run in parallel. Collectively, these three optimizers are referred to as \colorDistAll{\DistAll}.}
		\label{fig:CPUtime}
	\end{figure}

	We provide the CPU times required to solve the DO problems under consideration in Figure~\ref{fig:CPUtime}. A comparison of the solution times for \colorCentFinal{\CentFinal} and \colorCentInter{\CentInter} reveals that the addition of intermediate buffer tanks and implementation of a centralized monolithic optimizer results in a \SI{108}{\percent} increase in solution time. On the other hand, the total solution time for \colorDistAll{\DistAll} is \SI{85}{\percent} and \SI{93}{\percent} less than that of \colorCentFinal{\CentFinal} and \colorCentInter{\CentInter}, respectively, highlighting the effectiveness of implementing distributed optimizers and solving DO problems with smaller-scale DAE systems in terms of computational cost. Notice that \colorDist{\DistRSR} must be solved before \colorDist{\DistB} and \colorDist{\DistG}, which are solved in parallel afterward.
    	
	Furthermore, implementing a centralized monolithic optimizer for the entire process leads to DO problems with large-scale DAE systems. NLP convergence for such problems, particularly for \colorCentInter{\CentInter}, is highly sensitive to the initial guess and variable scaling, resulting in ill-conditioning issues and consequently non-convergence. However, for the \colorDistAll{\DistAll} distributed optimizers, the NLP solvers converge robustly and are easier to initialize and scale due to the significantly smaller DAE systems involved in their DO problems.
	
	Therefore, although the \colorDistAll{\DistAll} yields \SI{1}{\percent} less addition in the total operating profit than \colorCentInter{\CentInter}, it is preferable to implement the former, especially for online applications, due to the significant reduction in solution time and the easier convergence of the NLP solvers in its DO problems.

	\section{Conclusion} \label{sec:conclusion}

	We investigate the operational flexibility potential of a fully electrified biodiesel production process by proposing different process configurations using intermediate and final buffer tanks and solving offline DO problems. This process comprises reaction, separation, and recycle components, indicating the application potential of its operational flexibility for a broad range of chemical processes. By incorporating both intermediate and final buffer tanks, we fully exploit the flexibility potential of the process, leading to total energy savings of up to \SI{29}{\percent} relative to the SS benchmark. Compared to implementing only final buffer tanks, the use of both intermediate and final buffer tanks yields an additional \SI{8}{\percent} savings. Furthermore, the utilization of intermediate buffer tanks enables production rates to reach process equipment limits, thereby leveraging the full potential of power demand flexibility.
	
	In addition to the flexibility in production rates, imposing purity limits on the outlet streams of the buffer tanks allows for flexibility in the produced product purities. Distillation processes require varying power levels for different product grades. Thus, such flexible purity production enhances the overall power demand flexibility. We demonstrate this by producing glycerol with purities below $97.5$~wt\unit{\percent} while still complying with the required minimum $99$~wt\unit{\percent} limits on the delivered final product.
		
	We observe that the intermediate buffer tanks play a critical role in decoupling the dynamics of their upstream and downstream processes. We explore the significance of each buffer tank in providing flexibility and demonstrate how these tanks are filled or utilized based on the electricity price profile. Moreover, to ensure that the final state of the tanks is similar to the initial state, it is necessary to impose endpoint constraints on the liquid levels and content purities. This ensures that the optimizer does not exploit the available holdups and purities, thus enabling online application in a moving horizon fashion.
	
	Furthermore, we demonstrate that intermediate buffer tanks can be also used to facilitate the implementation of distributed optimizers for various process parts, rather than relying on a centralized monolithic optimizer for the entire process. This approach allows for solving DO problems with smaller DAE system sizes, resulting in improved computational performance. Our study shows that the distributed optimization approach leads to solution time savings of up to \SI{93}{\percent} compared to the centralized approach, while its NLP solvers are more robust and less prone to non-convergence. However, to apply the distributed optimization approach, it is necessary to tear the residual oil recycle stream and fix its flow rate, resulting in less efficient utilization of the oil feed. Nonetheless, we demonstrate that this loss in economic performance is small. Therefore, given the economic and computational performances of all considered DO cases, the distributed optimization approach is more favorable, particularly for online applications.
	
	Although our current process design lacks consideration of heat integration, it is important to note that the integration of heat across multiple units is common in complex and modern chemical plants. Therefore, future research should also investigate the implications of heat integration on process flexibility, particularly by exploring the degrees of freedom available to optimization. Additionally, to electrify the process, we utilize separate heat pumps for refrigeration and heating in reboilers, multiple heat pump stages, internal heating exchange, and multistage compression, based on literature research. In light of this, future work should focus on implementing process design optimization to determine optimal structures and configurations, as well as heat pump integration.

    We perform the sizing of process units and buffer tanks heuristically and based on literature findings. Instead one could formulate these also as optimization variables/problems to be solved alongside the DO problems for flexible operation. This would involve employing stochastic programming and presents an interesting avenue for future work. Nevertheless, it is important to note that addressing this aspect may entail dealing with significantly large problem sizes, given the requirement to consider multiple operational scenarios, particularly for electricity price profiles.
	
	As further future work, applying online DO for demand-side management, specifically economic nonlinear model predictive control, to the considered process would greatly enhance our understanding of the real-time applicability of the proposed flexibility-oriented process design and flexibilization strategies. The use of distributed optimizers promises good computational performance and should be considered for online applications. Additionally, hierarchical control structures with optimal scheduling and lower-level tracking control can be also explored as future work. Finally, in cases where optimization problems are computationally expensive and encounter convergence issues, it would be worthwhile to investigate model reduction techniques, such as incorporating surrogate models, and explore algorithm and implementation improvements.

	\section{Data availability}
	The Modelica process model is available open-source under \href{http://permalink.avt.rwth-aachen.de/?id=135903}{permalink.avt.rwth-aachen.de/?id=135903}.

	\section{Author contributions}
	
	\textbf{Mohammad El Wajeh}: Conceptualization, Methodology, Software, Investigation, Validation, Formal analysis, Writing - Original draft \\
	\textbf{Adel Mhamdi}: Conceptualization, Project administration, Supervision, Funding acquisition, Writing - Review and editing \\
	\textbf{Alexander Mitsos}: Conceptualization, Project administration, Supervision, Resources, Writing - Review and editing

    \section{Competing interests} 
	None.
	
	\section{Acknowledgments}
	The authors gratefully acknowledge the financial support of the Kopernikus project SynErgie by the German Federal Ministry of Education and Research (BMBF) and the project supervision by the project management organization Projektträger Jülich (PtJ). The authors thank Johannes Faust, Jan Schulze, Aron Zingler, and Eike Cramer from AVT.SVT, RWTH Aachen for fruitful discussions.

	\printbibliography[heading=bibnumbered]
	
\end{document}